\documentclass[reqno]{amsart}
\usepackage{amssymb,amsthm,amsfonts,amstext}
\usepackage{amsmath}

\makeatletter \@addtoreset{equation}{section} \makeatother

\renewcommand\thefigure{\thesection.\@arabic\c@figure}
\renewcommand\thetable{\thesection.\@arabic\c@table}
\newtheorem{theorem}{Theorem}[section]
\newtheorem{lemma}[theorem]{Lemma}
\newtheorem{proposition}[theorem]{Proposition}

\newtheorem{definition}[theorem]{Definition}
\newtheorem{remark}[theorem]{Remark}

 \theoremstyle{remark}
 
\newcommand{\mc}[1]{{\mathcal #1}}

\newcommand{\bb}[1]{{\mathbb #1}}

\newcommand{\<}{\langle}
\renewcommand{\>}{\rangle}

\title[Scaling limits for gradient systems in random environment]{SCALING LIMITS FOR GRADIENT SYSTEMS IN RANDOM ENVIRONMENT}
\date{\today}
\begin{document}
\author{Patr\'{\i}cia Gon\c{c}alves and Milton Jara}
\begin{abstract}
It is well known that the hydrodynamic limit of an interacting particle system satisfying a gradient condition (such as the zero-range process
or the symmetric simple exclusion process) is given by a possibly non-linear parabolic equation and the equilibrium fluctuations from this limit are
given by a generalized Ornstein-Uhlenbeck process.

We prove that in the presence of a symmetric random environment, these scaling limits also hold for almost every choice of the random
environment, with an homogenized diffusion coefficient that does not depend on the realization of the random environment.
\end{abstract}
\renewcommand{\subjclassname}{\textup{2000} Mathematics Subject Classification}
\begin{thanks} {The first author wants to thank
    F.C.T. (Portugal) for supporting her Phd with the grant /SFRH/ BD/ 11406/ 2002.}
\end{thanks}
\address{IMPA, Estrada Dona Castorina 110, CEP 22460-320, Rio de Janeiro, Brazil}
\email{patg@impa.br, monets@impa.br}
 \maketitle

\section{Introduction}

Consider a system of particles evolving on a multidimensional, periodic integer lattice of period $2N$. Each particle performs a continuous-time random
walk with rates $p(x,y)$ that depend on both the position $x$ and the destination site $y$. These rates are
chosen as a fixed realization of a random field, in such a way that the resulting single-particle random walk is
reversible with respect to the counting measure on the lattice. We call these rates the {\em random environment}.

Particles interact between them only when they share a site,
through an interaction function $g: \bb N_0 \to \bb R_+$. The dynamics for this system is the following. At each time $t$,
let $\eta_t(x)$ denote the number of particles at the site $x$. For each pair of sites $\<x,y\>$, after an exponential waiting time of
rate $g(\eta_t(x)) p(x,y)$ the particle at site $x$ jumps to site $y$. This is done independently for each pair $\<x,y\>$ and
after each jump, the exponential waiting time for each pair $\<x,y\>$ starts afresh.

Such a system can be understood as a model for diffusion in heterogeneous media. The purpose of this article is to
study the scaling limits of this system as $N \to \infty$ and mostly the influence
of the randomness in this limit. As we will see, when the underlying random field is ergodic, stationary and satisfies an ellipticity condition, for
any realization of the random environment the scaling limit depends on the randomness only through some constants
which depend on the distribution of the random transition rates, but not on the particular realization of the random environment.

In this article we study two related scaling limits for this process: the hydrodynamic limit and the equilibrium fluctuations. The first one is
a law of large numbers for the empirical distribution of particles when the process starts from a configuration of particles with macroscopic
density close to some initial profile while the second one is a central limit theorem for the empirical distribution of particles
 when the system starts from an equilibrium measure.

The hydrodynamic limit has been obtained in the context of exclusion processes in \cite{FM} when the dimension $d \geq 3$ and extended in
\cite{Qua} to any dimension. In these references, it is not assumed the reversibility for the one-particle random walk with respect to the
counting measure on the lattice, so in this sense their results are more general than ours. Their approach is based on the generalization of the
non-gradient method of Varadhan \cite{Var}, \cite{Qua2} for the case of random transition rates. In \cite{Fag}, \cite{JL} the one-dimensional
simple exclusion process is considered.

In the reversible situation, we introduce the {\em corrected empirical process}. This process satisfies the gradient condition, which is a key
property from which hydrodynamics and equilibrium fluctuations can be easily obtained like in the non-random situation \cite{GPV}, \cite{Cha}.
Therefore, our approach is simpler, does not require any mixing condition and can be generalized to situations in which the non-gradient method does not
apply, like kinetically constrained particle systems, the zero-range process with bounded interaction rate and particle systems in
non-homogeneous lattices \cite{Jar}.

The introduction of a corrected empirical measure can be understood as a version of Tartar's  compensated compactness lemma in the context of particle
systems. In this reversible situation the averaging due to the dynamics and the inhomogeneities introduced by the random media factorize after
introducing the corrected empirical process, in such a way that we can average them separatedly. For the dynamic averaging, we use the entropy
method of \cite{GPV} to derive the hydrodynamic limit, while for the equilibrium fluctuations we adopt Chang's proof of the Boltzmann-Gibbs
principle \cite{Cha}; for the averaging of the random environment we use $\Gamma$-convergence. With this procedure, the scaling limits of the
corrected empirical process are obtained. After this, we prove that in the limit as $N \to \infty$, the corrected empirical process and the
original empirical process are close enough to recover the scaling limit for the original empirical distribution of particles.

In order to see how far can this picture be taken, we also prove the Boltzmann-Gibbs
principle for functions that depend on both the particle configuration and the random environment. Notice
that this more general version of the Boltzmann-Gibbs principle is not needed to obtain the equilibrium fluctuations
for the empirical density of particles.

The Boltzmann-Gibbs principle states that non-conserved quantities oscillate
faster than conserved quantities, and therefore when averaged in time, only the projections over the density field are
observed. In consequence, the Boltzmann-Gibbs principle is interesting by its own. In order to give further motivations for the study
of the Boltzmann-Gibbs principle for random functions, we present two applications at the end of the article.

The article is structured as follows. In Section \ref{s1} we describe the model and the main results. Section \ref{s2} is devoted to the
proof of the hydrodynamic limit for this process and in the subsequent section we present the equilibrium fluctuations. The proof of the
Boltzmann-Gibbs principle is referred to Section \ref{s4}. For the reader's convenience, we include some well-known, but rather technical lemmas and definitions in
the Appendix.

\section{Notations and results}
\label{s1}
\subsection{The zero-range process} 
We define the zero-range process as a continuous-time Markov process $\eta_t$ with state space $\Omega_{N}^{d} = \{\eta: \bb T^d_N \to \bb N_0\}$, where $\bb T_{N}^{d}$
is the $d$-dimensional
discrete torus $N^{-1} \bb Z^d/ 2N\bb Z^d$. We consider $\bb T_{N}^{d}$ as a subset of $U^{d} = [-1,1]^d$ with periodic boundary conditions. This
process has a generator whose action over local functions  $f: \Omega_N^{d} \to \bb R$  is given by
\[
\mc L_N f(\eta) = \sum_{x, y \in \bb T_N} p_N(x,y) g\big(\eta(x)\big) \big[f(\eta^{xy})-f(\eta)\big],
\]
where $p_N: \bb T_N^{d} \times \bb T_N^{d} \to \bb R_+$ is the jump rate of a random walk in $\bb T_N^d$, $g: \bb N_0 \to \bb R_+$ is the
interaction rate between the particles and $\eta^{xy} \in \Omega_N^d$ is given by
\[
\eta^{xy}(z) =\begin{cases}
\eta(x)-1, z =x \\
\eta(y)+1, z=y\\
\eta(z) , z \neq x,y\\
\end{cases}.
\]

Notice that the dynamics of $\eta_t$ conserves the number of particles. In particular, the process $\eta_t$ is well defined for any initial
configuration $\eta_0 \in \Omega_N^d$, since in that case the state space is finite.

We will assume that the interaction rate $g$ has linear growth:
\begin{equation}\label{lineargrowth}
 \exists c_0>0: c_0^{-1}n \leq g(n) \leq c_0
n \quad\ \forall n \in \bb N_0
\end{equation}

We will also assume that the motion of a single particle is a nearest-neighbor random walk, so we take $p_N(x,y)=0$ if $|x-y| \neq 1/N$,
where $|x-y| = \sum_{i \leq d} |x_i -y_i|$ is the sum norm in $\bb R^d$. This last hypothesis is not essential, but it simplifies the notation.
We further assume that $p_N(x,y) =p_N(y,x)$ for all $x, y \in \bb T_N^d$. This hypothesis will ensure the reversibility of the process $\eta_t$
with respect to the measures $\nu_\rho$ defined below, and the reversibility of the randomm walk generated by $p_N(x,y)$, which is crucial in what follows.

For each $\alpha \geq 0$, let $\bar \nu_\alpha$ be the product measure in $\Omega_N^d$ whose marginals are given by
\[
\bar \nu_\alpha\big\{\eta;\eta(x)=k\big\} =\frac{1}{Z(\alpha)}
\frac{\alpha^k}{g(k)!},
\]
where $g(k)!=g(1)\cdots g(k)$ for $k\geq 1$, $g(0)=1$ and $Z(\alpha)$ is the normalizing constant for which $\bar \nu_\alpha(\Omega_N^d)=1$. By
the linear growth of $g$ (\ref{lineargrowth}), $\bar \nu_\alpha$ is well defined for all $\alpha\geq 0$.

Define $\rho=\rho(\alpha)$ as the density of particles with respect to $\bar \nu_\alpha$, namely:
\[
\rho(\alpha) = E_{\bar \nu_\alpha}\big[ \eta(x) \big] = \sum_{k \geq 0}
\frac{1}{Z(\alpha)} \frac{k \alpha^k}{g(k)!} = \frac{\alpha
Z'(\alpha)}{Z(\alpha)}.
\]

Again by the linear growth of $g$, $\alpha \mapsto \rho(\alpha)$ is an homeomorphism from $[0,\infty)$ to $[0,\infty)$ and the inverse function
$\alpha= \alpha(\rho)$ is well defined for all $\rho \in [0,\infty)$. We define $\nu_\rho =\bar \nu_{\alpha(\rho)}$ and
$\phi(\rho)=E_{\nu_{\rho}}[g(\eta(0))]$. Due to the symmetry of $p_N(x,y)$, the measure $\nu_\rho$ is invariant and reversible for this process.

\subsection{The random environment}
Now we discuss the choice of the jump rates $p_N(x,y)$. Let $(\mc X, \mc F, P)$
be a probability space and take a family $\{\theta_x; x \in \bb Z^d\}$ of $\mc
F$-measurable mappings $\theta_x: \mc X \to \mc X$ such that
\begin{itemize}
\item[i)] $P(\theta_x^{-1}A) = P(A)$ for all $A \in \mc F$, $x \in \bb Z^d$.
\item[ii)] $\theta_z \theta_{z'} = \theta_{z+z'}$ for all $z, z' \in \bb Z^d$.
\item[iii)] If $\theta_z A = A$ for all $z \in \bb Z^d$, then $P(A) = 0$ or $1$.
\end{itemize}

In this case we say that the family $\{\theta_x\}_{x\in{\bb{Z}^d}}$ is invariant and ergodic under $P$. Let $a= (a_1,...,a_d): \mc X \to \bb
R^d$ be a $\mc F$-measurable function such that there exists $\epsilon_0>0$ with
\begin{equation} \label{boundconditiononRE}
\epsilon_0 \leq a_i(\omega) \leq \epsilon_0^{-1} \text{ for all }\omega \in \mc X \text{ and } i=1,...,d.
\end{equation}

 Fix $\omega \in \mc X$. For each $x \in \{-1+1/N,-1+2/N,...,1\}^d$ and $i=1,...,d$, define
\begin{equation}\label{RE}
p_N(x,x+e_i/N)= p_N(x+e_i/N, x)= N^2 a_i(\theta_{Nx}\omega)
\end{equation}
to which we call the {\em random environment}.

For each $G: \bb T_N^d \to \bb R$, define the operator $L_N G$ by
\[
L_N G(x) = \sum_{y \in \bb T_N^d} p_N(x,y) \big[G(y) -G(x)\big].
\]

In the space of functions $l_N(\bb T_N^d)=\{f: \bb T_N^d \to \bb R\}$, define the following norms:
\[
||f||_{0,N}^2 = \frac{1}{N^d} \sum_{x \in \bb T_N^d} f(x)^2
\]
and
\[
||f||_{1,N}^2 = ||f||^{2}_{0,N} + \frac{1}{N^d} \sum_{\substack{x, y \in \bb T_N^d\\|x-y| =1/N}} N^2 \big[f(y) -f(x)\big]^2.
\]

We denote by $\mc L^2_N$ the space of functions $l_N(\bb T_N^d)$ endowed with the norm $||\cdot||_{0,N}$ and by $\<\cdot,\cdot\>_N$ the inner
product in $\mc L^2_N$. Define $\mc H_{1,N}$ as the space of functions in $l_N(\bb T_N^d)$ endowed with the norm $||\cdot||_{1,N}$.

Denote by $\mc L^2(U^d)$ the space of square integrable functions in $U^d$ with respect to the Lebesgue measure and by $||\cdot||_0$ the corresponding norm in
$\mc L^2(U^d)$. For each $k \geq 1$, denote by $\mc H_k(U^d)$ the Sobolev space in $U^d$ defined as the completion of $\mc C_\infty(U^d)$ under
the norm
\[
||f||_k^2 = \sum_{|\alpha| \leq k} ||\partial^\alpha f||_0^2,
\]
where $|\alpha|$ denotes the order of the multi-index $\alpha$ and
$\partial^\alpha$ is the partial derivative of order $\alpha$.

The definition of convergence of a sequence $f_N \in \mc H_{1,N}$ (or $\mc L_N^2$) to $f \in \mc H_1(U^d)$ (or $\mc L^2(U^d)$) is given in
Appendix \ref{ap1}.

From the homogenization theory the following holds:

\begin{proposition} \label{p1}
Fix a typical realization of $p_N(\cdot,\cdot)$ and $\lambda >0$. There exists a positive defined matrix $\mc A$ that depends only on the
distribution of $a=(a_1,...,a_d)$ such that for any $f_N$ and $f$ such that $f_N \in \mc H_{-1,N}$ converges strongly to $f \in \mc
H_{-1}(U^d)$, $u_N$ converges weakly in $\mc H_{1,N}$ to $u$, where $u_N$ is defined as the solution of the equation
\[
\lambda u_N(x) - L_N u_N(x) = f_N(x)
\]
and $u$ is the solution of the equation
\[
\lambda u(x) - \nabla \cdot \mc A \nabla u(x) = f(x).
\]
\end{proposition}

A proof of this proposition can be found in \cite{pr}. Notice that the statement of this proposition makes sense for any choice of the jump rate
$p_N(x,y)$.

In order to prove the hydrodynamic limit we need this property on the jump rates $p_N(x,y)$ and for this reason we introduce the following
definition.

\begin{definition}
We say that a family of jump rates $\{p_N: \bb T_N^d \times \bb T_N^d \to \bb R_+\}_N$ admits homogenization, if there exist a constant
$\epsilon_0
>0$ such that $\epsilon_0 N^2 \leq p_N(x,y) \leq \epsilon_0^{-1} N^2$ and a
matrix $\mc A$ such that for any $f \in \mc H_{-1}(U^d)$ smooth enough there exists a sequence $f_N$ converging strongly in $\mc H_{-1,N}$ to
$f$ such that the solution $u_N \in \mc H_{1,N}$ of the equation
\[
\lambda u_N(x) - L_N u_N(x) = f_N(x)
\]
converges weakly in $\mc H_{1,N}$ to the solution $u \in \mc H_1(U^d)$ of
\[
\lambda u(x) - \nabla \cdot \mc A \nabla u(x) = f(x).
\]

In this case, we say that the matrix $\mc A$ is the $\Gamma$-limit of $L_N$.
\end{definition}

For our purposes, $f$ will be smooth enough if it is three times countinuously differentiable.

\begin{remark}
By the theory of $\Gamma$-convergence, the matrix $\mc A$ satisfies the coerciveness assumption $\epsilon_0 |\xi|^2 \leq \sum_{ij} \xi_i \xi_j
\mc A_{ij} \leq \epsilon_0^{-1} |\xi|^2$ for all vectors $\xi \in \bb R^d$. In the previous definition, nothing excludes the possibility of the
matrix $\mc A$ to be a function of the position $x \in U^d$. See \cite{os} for a one-dimensional example on which the $\Gamma$-limit of $L_N$ is
not constant in space.
\end{remark}

\subsection{Hydrodynamic limit}
Fix a function $\rho_0:U^d \to \bb R_+$. A family of measures $\{\mu_N\}_{N\geq{1}}$ in $\Omega_N^d$ is said to be associated to the profile $\rho_0$ if
for any function $G \in \mc C(U^d)$ and any $\epsilon
>0$,
\[
\lim_{N \to \infty} \mu_N \Big(\eta \in \Omega_N^d; \Big|\frac{1}{N^{d}} \sum_{x \in \bb T_N^d} \eta(x) G(x) - \int \rho_0(x) G(x) dx \Big| >
\epsilon \Big) =0.
\]

Here and in the sequel, denote by $E_{\mu_N}$ the expectation with respect to $\mu_N$ and by $\bb E_{\mu_N}$ the expectation with respect to
$\bb P_{\mu_N}$, the distribution of the process $\eta_t$ starting from $\mu_N$ in $\mc D([0,T], \Omega_N^d)$. We follow the evolution of the
process $\eta_t$ in a finite time interval $[0,T]$ in order to avoid uninteresting complications due to the lack of compactness of $[0,\infty)$.

Let $\rho >0$ be a fixed density.
The entropy of $\mu_N$ with respect to $\nu_\rho$ is defined by

\[
H_N(\mu_N|\nu_\rho) =
\begin{cases}
\int \frac{d\mu_N}{d\nu_\rho} \log \frac{d\mu_N}{d\nu_\rho} d \nu_\rho, \text{
if } \mu_N \ll \nu_\rho\\
+ \infty, \text{ otherwise}
\end{cases},
\]
where for two measures $\mu$ and $\nu$, $\mu\ll\nu$ means that the measure $\nu$ is absolutely continuous with respect to $\mu$.

We introduce a partial order $\preceq$ in $\Omega_N^d$ as follows. For $\eta$, $\eta'$ in $\Omega_N^d$, we say that $\eta \preceq \eta'$ if
$\eta(x) \leq \eta'(x)$ for every $x \in \bb T_N^d$. Once there is a partial order in the space state $\Omega_{N}^{d}$, we can introduce a
partial order in the space of measures in $\Omega_{N}^{d}$.  We say that $\mu_N$ is stochastically dominated by $\nu_\rho$ (also denoted by
$\mu_N \preceq \nu_\rho$) if there exists a measure $\bar \mu$ in $\Omega_N^d\times\Omega_N^d$ such that:

\quad

\begin{itemize}
\item[i)] For all $\eta \in \Omega_N$, $\bar \mu(\eta,\Omega_N) = \mu_N(\eta)$.
\item[ii)] For all $\eta \in \Omega_N$, $\bar \mu(\Omega_N, \eta) =
\nu_\rho(\eta)$.
\item[iii)] The set $\{(\eta,\eta'); \eta \preceq \eta'\}$ has full measure
under $\bar \mu$.
\end{itemize}

 In this case we say that $\bar \mu$ is a {\em coupling} of $\mu_N$ and
 $\nu_{\bar \rho}$.

\begin{theorem}
Let $\rho_0:U^d \to \bb R$ be a bounded profile, and let $\{\mu_N\}_{N\geq{1}}$ be a sequence of measures in $\Omega_N^d$ associated to the
profile $\rho_0$. Assume that the interaction rate $g(\cdot)$ is non-decreasing and has linear growth (see Section \ref{lineargrowth}). Suppose that
there exist constants $K_0$ and $\bar \rho$ such that $H(\mu_N|\nu_{\bar \rho}) \leq K_0N^d$ and $\mu_N \preceq \nu_{\bar \rho}$ for every $N$ large enough. Suppose also
that the jump rates $p_N(x,y)$ admit homogenization with homogenized matrix $\mc A$.

Then, for every $t \leq T$, every continuous function $G : U^d \to \bb R$ and every $\delta
>0$,
\[
\lim_{N \to \infty} \bb P_{\mu_N} \Big[ \big| \frac{1}{N^{d}} \sum_{x \in \bb T_N^d} G(x) \eta_t(x) - \int G(u) \rho(t,u) du \big| > \delta
\Big]=0,
\]
where $\rho(t,u)$ is the unique weak solution of the hydrodynamic equation
\begin{equation}
\begin{cases}
\partial_t \rho = \nabla \cdot \big( \mc A \nabla \phi(\rho) \big) \\
\rho(0,\cdot) = \rho_0(\cdot).
\end{cases}
\label{echid}
\end{equation}
\label{t1}
\end{theorem}

In the sake of completeness we introduce the definition of weak solutions of equation (\ref{echid}).
\begin{definition}
Fix a bounded profile $\rho_{0}:U^{d}\rightarrow{\mathbb{R}}$. A bounded function $\rho:[0,T]\times U^{d}\rightarrow{\mathbb{R}}$ is a
{\em weak solution} of equation (\ref{echid}) if for every function $G:[0,T]\times U^{d}\rightarrow{\mathbb{R}}$ of class
$C^{1,2}([0,T]\times U^{d})$,
\begin{equation*}
\int_0^t \int_{U^d} \Big\{ \rho(s,u) \partial_s G(s,u) + \phi(\rho(s,u))\nabla \cdot \mc A \nabla G(s,u) \Big\}du ds
\end{equation*}
\begin{equation}
+\int_{U^d}\rho_{0}(u)G(0,u)du=\int_{U^d}\rho(T,u)G(T,u)du. \label{definitionweaksolution}
\end{equation}
\end{definition}

Let $\mc M_+$ be the set of positive Radon measures in $U^d$. The empirical
measure $\pi_t^N$ is defined as the process in $\mc D([0,T],\mc M_+)$  given by

\[
\pi_t^N(du)=\frac{1}{N^d} \sum_{x \in \mathbb{T}_N^d} \eta_t(x) \delta_x(du),
\]
where $\delta_x$ is the Dirac distribution at $x$.

For $G: U^d \to \bb R$ continuous, define $\pi_t^N(G) = \int G(u) \pi_t^N(du)$. The statement of Theorem \ref{t1} is equivalent to say that
under $\bb P_{\mu_N}$ the random variables $\pi_t^N(G)$ converge in probability to $\int G(u) \rho(t,u) du$ for every $G$ continuous and every
$t \in [0,T]$. We will prove a stronger result for $\pi_t^N$:

\begin{theorem}
\label{t2}
Under the hypothesis of Theorem \ref{t1}, $\pi_t^N$ converges in distribution in
$\mc D([0,T],\mc M_+)$ to the trajectory $\rho(t,u) du$.
\end{theorem}

\begin{remark}
\label{r1}
Since $\rho(t,u) du$ is a deterministic element of $\mc D([0,T],\mc M_+)$, the
convergence in distribution of $\pi_t^N$ implies its convergence in probability,
from which Theorem \ref{t1} follows.
\end{remark}

\subsection{Equilibrium fluctuations}

Now we state a central limit theorem for the empirical measure, starting from an equilibrium measure $\nu_\rho$. Fix $\rho>0$ and denote by
$\mc S(U^d)$ the Schwartz space of infinitely differentiable functions in $U^d$.

Denote by $\mathcal{Y}_{\cdot}^{N}$ the density fluctuation field, a linear functional acting on functions $G \in \mc S(U^d)$ as
\begin{equation}\label{densityfield}
\mc Y_t^N(G) = \frac{1}{N^{d/2}} \sum_{x \in \bb T_N^d} G(x) (\eta_t(x) -\rho).
\end{equation}

Notice that
\begin{equation*}
\mc Y_t^N(G) = N^{d/2}\Big(\int G(u) \pi_t^N(du)-\rho \int G(u)du\Big).
\end{equation*}
In this way we have defined a process in $\mc D([0,T], \mc S'(U^d))$, where $\mc S'(U^d)$ is the space of tempered distributions, which
corresponds to the dual of the Schwartz space $\mc S(U^d)$.

\begin{theorem}
\label{tfeq} Consider the fluctuation field $\mc Y_{\cdot}^N$ defined above. Assume that the interaction rate $g(\cdot)$ has linear growth and
that the jump rates admit homogenization with homogenized matrix $\mc A$.

Then, for every $t_1,...,t_k \in [0,T]$ and every $G_1,...,G_k \in \mc S(U^d)$, the vector $(\mc Y_{t_1}^N(G_1),...,\mc Y_{t_k}^N(G_k))$
converges in distribution to $(\mc Y_{t_1}(G_1),..., \mc Y_{t_k}(G_k))$, where $\mc Y_t$ is the generalized Ornstein-Uhlenbeck process of
characteristics $\phi'(\rho)\nabla \cdot \mc A\nabla$ and $\sqrt{\phi(\rho) \mc A} \nabla$.
\end{theorem}

\section{Proof of Theorem \ref{t1}}
\label{s2} By remark \ref{r1}, in order to prove Theorem \ref{t1} it is enough to prove Theorem \ref{t2}. The proof of Theorem \ref{t2} follows
the standard lines of the proof of hydrodynamic limit by the entropy method for interacting particle systems. The route to proceed is the
following:

First we show that the distributions of $\pi_t^N$ in $\mc D([0,T],\mc M_+)$ form a tight sequence. Then we prove that the limit points of
$\pi_t^N$ are concentrated on trajectories of measures absolutely continuous with respect to the Lebesgue measure in $U^d$ with a bounded
density. Finally, we prove that these limit points are concentrated on weak solutions of the hydrodynamic equation (\ref{echid}). By the
uniqueness of these weak solutions on the space of bounded functions we conclude that $\pi_t^N$ has a unique limit point, concentrated on the
trajectory with density $\rho(t,u)$, where $\rho(t,u)$ is the weak solution of equation (\ref{echid}). Since the topology of convergence in
distribution is metrizable, we conclude that the whole sequence $\pi_t^N$ converges to $\rho(t,u)du$.

Unfortunately, this plan cannot be accomplished directly for $\pi_t^N$, but for another auxiliary process, the {\em corrected empirical
measure}, that we define below.

Let $\lambda>0$ be fixed. A function $G : U^d \to \bb R$ is said to be regular if the function $f_N \in \mc L_N^2$ defined by $f_N(x)
= \lambda G(x) - \nabla \cdot \mc A \nabla G(x)$ converges strongly in $\mc H_{-1,N}$ to $\lambda G - \nabla\cdot \mc A \nabla$. Notice that a
sufficient condition for $G$ to be regular, is $G \in \mc C^3(U^d)$, where $\mc C^3(U^d)$ denotes the space of three times continuously
differentiable functions on $U^d$.

Let $G: U^d \to \bb R$ be regular. For each $N \geq 1$, define $R^\lambda G(x) = \lambda G(x) - \nabla \cdot \mc A \nabla G(x)$ and
$G_N^\lambda: \bb T_N^d \to \bb R$ as the solution of

\begin{equation} \label{poisson1}
\lambda G_N^\lambda(x) - L_N G_N^\lambda(x) = R^\lambda G(x).
\end{equation}
By Lemma \ref{l2}, the following estimates hold:
\begin{equation} \label{esti1}
||G^{\lambda}_{N}||_{0,N}\leq{\lambda^{-1}||R^{\lambda}G||_{0,N}},
\end{equation}
\begin{equation} \label{esti2}
\frac{1}{N^{d}} \sum_{x,y\in{\mathbb{T}_{N}}} p_{N}(x,y) [G^{\lambda}_{N}(y)-G^{\lambda}_{N}(x)]^{2} \leq \lambda^{-1}||R^{\lambda}G||^{2}_{0,N}
\end{equation}
and
\begin{equation}
\label{esti3}
||G_N^\lambda||_{\infty, N} \leq \lambda^{-1} ||R^\lambda G||_{\infty,N}.
\end{equation}

We define the corrected empirical measure $\pi_t^{N,\lambda}$ by
\[
\pi_t^{N,\lambda}(G) = \frac{1}{N^d} \sum_{x \in \bb T_N^d} \eta_t(x) G_N^\lambda(x).
\]

Notice that $\pi_t^{N,\lambda}(G)$ is defined only for $G$ regular, so
$\pi_t^{N,\lambda}$ is not a well defined process
 in $\mc D([0,T],\mc M_+)$. Lemma \ref{l2} shows that $\pi_t^{N,\lambda}$ is a
 well defined process in the Sobolev space $\mc H^{-k}(U)$ for $k \geq 3$.
 However, this point will not be relevant for our proof of Theorem \ref{t2}.

Since $\mc M_+$ is separable and the vague topology in $\mc M_+$ is metrizable, in order to prove tightness of $\pi_t^N$ in $\mc D([0,T],\mc
M_+)$, it is enough to show tightness of $\pi_t^N(G)$ in $\mc D([0,T], \bb R)$ for $G$ in a dense subset
 of the set $\mc C(U^d)$ of continuous functions in $U^d$. Therefore, it is enough
 to prove tightness of $\pi_t^N(G)$ for $G$ regular.

 By Dynkin's formula,
\begin{equation}
\mc M_t^N(G) = \pi_t^{N,\lambda}(G) - \pi_0^{N,\lambda}(G) - \int_0^t\frac{1}{N^d} \sum_{x \in \bb T_N^d} g\big(\eta_s(x)\big) L_N
G_N^\lambda(x) ds \label{ec2}
\end{equation}
is a martingale of quadratic variation given by
\[
\<\mc M_t^N(G)\> = \int_0^t \frac{1}{N^{2d}} \sum_{x,y \in \bb T_N^d} g\big(\eta_s(x)\big) p_N(x,y) \big[G_N^\lambda(y) -
G_N^\lambda(x)\big]^2ds.
\]

We claim that $\mc M_t^N(G)$ goes to 0 as $N \to \infty$ in $\mc L^2(\bb
P_{\mu_N})$. In fact,
\[
\begin{split}
\bb E_{\mu_N}\big[\mc M_t^N(G)^2\big] &= \bb E_{\mu_N}\big[\<\mc M_t^N(G)\>\big]
\\
&= \int_0^t \frac{1}{N^{2d}} \sum_{x,y \in \bb T_N^d} \bb E_{\mu_N} \big[g(\eta_s(x))\big] p_N(x,y) \big[ G_N^\lambda(y) -G_N^\lambda(x)\big]^2
ds
\\
&\leq \frac{t \phi(\bar \rho)}{N^d} \lambda^{-1} ||R^\lambda G||_{0,N}^2
\xrightarrow{N \to \infty} 0.
\end{split}
\]
In order to obtain this last bound, we have used the estimate (\ref{esti2}), the
fact that $\mu_{N}$ is stochastically dominated by
$\nu_{\bar{\rho}}$ and Proposition \ref{p2}.

To prove tightness for the martingale $\mc M_t^N(G)$, we use the following
criterion, due to Aldous:

\begin{proposition}
A sequence of probability measures $\{P_N\}_N$ in $\mathcal
D([0,T],\mathbb R)$ is tight if
\begin{itemize}
\item[(i)] For all $0 \leq t\leq T$ and for all $\epsilon > 0$ there exists a
  finite constant $A$ such that  $\sup_N P_N(|x_t|>A) < \epsilon$,
\item [(ii)] For all $\epsilon > 0$,
\begin{equation*}
\lim_{\delta \to 0} \limsup_{N \to \infty}
\sup_{\substack{\tau \in \mathcal T\\\beta \leq \delta}}
P_N(|x_{\tau + \beta}-x_\tau|>\delta)=0\;,
\end{equation*}
where $\mathcal T$ is the set of stopping times with respect to the
canonical filtration, bounded by $T$.
\end{itemize}
\end{proposition}

A proof of this lemma can be found in \cite{kl}. In our case, condition $i)$ follows from the fact that $\mc M_t^N(G)$
 converges to $0$ in $\mc L^2(\bb P_{\mu_N})$ and Tchebyshev's inequality. On the other hand, by Doob's optimal sampling theorem, we have that
\[
\begin{split}
\bb P_{\mu_N} \big[ \big|
    &\mc M_{\tau+\beta}^N(G) -\mc M_{\tau}^N(G)\big| > \epsilon \big] \leq \\
    & \leq \frac{1}{\epsilon^2}\bb E_{\mu_N} \Big[ \<\mc M_{\tau+\beta}^N(G)\> - \<\mc M_{\tau}^N(G)\>\Big] \ \\
    &\leq \frac{1}{\epsilon^2} \bb E_{\mu_N}\Big[ \int_{\tau}^{\tau+\beta} \frac{1}{N^{2d}} \sum_{x,y \in \bb T_N^d} g\big(\eta_s(x)\big)p_N(x,y)  \big[G_N^\lambda(y)-G_N^\lambda(x)\big]^2 ds \Big] \\
    &\leq \frac{\beta C(G,c_0,\lambda,\epsilon_0)}{N^{d-2}} E_{\mu_N}\Big[ \frac{1}{N^d} \sum_{x \in \bb T_N^d} \eta(x) \Big].\\
\end{split}
\]
In this last bound we have used the conservation of the number of particles, the estimate (\ref{esti2}), and the uniform bound for $p_N(x,y)$.
Since the expected initial density of particles is bounded by $\bar \rho$, condition $ii)$ follows.

Notice that the integral term in (\ref{ec2}) can be written as
\begin{equation}
\int_0^t \frac{1}{N^d} \sum_{x \in \bb T_N^d} g\big(\eta_s(x)\big) \big[ \lambda G_N^\lambda(x) - R^\lambda G(x) \big] ds. \label{ec3}
\end{equation}

We see that
\[
\begin{split}
\bb E_{\mu_N} \Big[  \sup_{|s-t| \leq \delta}\Big|\int_{s}^{t} \frac{1}{N^d} \sum_{x \in \bb T_N^d} g\big(\eta_{t'}(x)\big) &\big[\lambda G_N^\lambda(x) - R^\lambda G(x)\big]dt'\Big|^2\Big] \leq  \\
&\leq \delta C(G,g) \bb E_{\mu_N} \big[ \sup_{t \in [0,T]} ||\eta_t||_{0,N}^2\big],
\end{split}
\]
that goes to 0 as $\delta \to 0$, uniformly in $N$ by Lemma \ref{l1}. Therefore, by Arzel\`a-Ascoli criterion,  the integral terms in (\ref{ec2}) form a tight sequence in $\mc D([0,T],\bb R)$ and their limit points are concentrated on continuous trajectories. By equation (\ref{ec2}) the sequence
$\pi_t^{N,\lambda}(G)$ is tight in $\mc D([0,T], \bb R)$. On the other hand, since $\mc M_t^N(G)$ goes to
 0 in $\mc L^2(\bb P_{\mu_N})$, any limit point of $\mc M_t^N(G)$ has null finite-dimensional distributions. Therefore, $\mc M_t^N(G)$
 converges to $0$ in distribution as a process in $\mc D([0,T],\bb R)$. Consequently, the limit points of
  $\pi_t^{N,\lambda}(G)$ are concentrated on continuous trajectories.

 Notice now that
\[
\bb E_{\mu_N} \big[ \sup_{t \in [0,T]} \big|\pi_t^{N,\lambda}(G) -\pi_t^N(G)\big|^2\big] \leq ||G_N^\lambda - G||_{0,N}^2 \bb E_{\mu_N} \big[ \sup_{t \in [0,T]} ||\eta_t||_{0,N}^2\big].
\]

By Proposition \ref{p1}, $||G_N^\lambda-G||_{0,N}$ converges to 0 as $N \to \infty$, and by Lemma \ref{l1} $\bb E_{\mu_N} \big[\sup_t
||\eta_t||_{0,N}^2\big]$ is bounded  in $N$.  Therefore, $\sup_t|\pi_t^{N,\lambda}(G)-\pi_t^N(G)|\to 0$ in $\mc L^2(\bb P_{\mu_N})$. A simple
$\varepsilon/3$ argument   allows us to obtain from this result that $\pi_t^N(G)$ is also tight in $\mc D([0,T],\bb R)$ and that $\pi_t^N(G)$
and $\pi_t^{N,\lambda}$ have the same limit points. Since the set of regular functions is dense in $\mc C(U^d)$, this ends the proof of
tightness for $\pi_t^N$ in $\mc D([0,T],\mathcal{M}_{+})$.

Let $\pi_t$ be a limit point of $\pi_t^N$, and let $Q$ be its distribution in $\mc D([0,T],\mc M_+)$. For any positive function $G \in \mc
C(U^d)$,
\[
\begin{split}
Q(\pi_t(G) > M) &\leq \liminf_{N \to \infty} Q(\pi_t^N(G) > M) \\
&= \liminf_{N \to \infty} \mu_{N} \big(N^{-d} \sum_{x \in \bb T_N^d} \eta_{t}(x)
G(x) > M\big) \\
&\leq \liminf_{N \to \infty} \nu_{\bar \rho} \big(N^{-d} \sum_{x \in \bb T_N^d}
\eta(x) G(x) > M\big) \\
&\leq \mathbf{1}{\Big( \int G(u) du > M/\bar \rho\Big)}.
\end{split}
\]

Here we have used once more, the fact of $\mu_{N}$ being stochastically dominated by an invariant measure $\nu_{\bar\rho}$ and Proposition
\ref{p2}.

Therefore, if $0 \leq G \leq 1$ then $Q(\pi_t(G) > 2^d \bar \rho)=0$. By the dominated convergence theorem, for every closed $B \subseteq U^d$
it holds that $Q(\pi_t(B) > \bar \rho \Lambda(B))=0$, where $\Lambda$ denotes the Lebesgue measure in $U^d$. In particular, the process $\pi_t$
is concentrated on measures absolutely continuous with respect to $\Lambda$.

Let $\pi(t,u)$ be the density of $\pi_t$ with respect to $\Lambda$. The same estimates prove that $\pi(t,u)$ is bounded by $\bar \rho$ in
$[0,T]\times U^d$.

Notice that $R^\lambda G$ is a smooth function, but $\lambda G_N^\lambda(x)$ it
is not smooth. However,
\[
\begin{split}
\bb E_{\nu_N}\Big[ \Big| \int_0^t \frac{1}{N^d} \sum_{x \in \bb T_N^d} g\big(\eta_s(x)\big) &\big[G_N^\lambda(x)-G(x)\big]ds\Big|^2\Big] \leq \\
&\leq c_0^{-2} t \int_0^t \bb E_{\mu_N}\Big[ \Big(\frac{1}{N^d} \sum_{x \in \bb T_N^d}
\eta_s(x)\big|G_N^\lambda(x) -G(x)\big|\Big)^2\Big]ds \\
&\leq c_0^{-2} t^2 \int \Big(\frac{1}{N^d} \sum_{x \in \bb T_N^d} \eta(x)
\big|G_N^\lambda(x)-G(x)\big|\Big)^2d\nu_{\bar \rho}(\eta) \\
&\leq c_0^{-2} t^2 \int \eta(0)^2 d\nu_{\bar \rho} ||G_N^\lambda -G||^2_{0,N}
\xrightarrow{N \to \infty} 0.
\end{split}
\]

In the previous we used Schwarz inequality together with the translation invariance of $\nu_{\bar\rho}$. As a consequence,
\begin{equation} \label{ec4}
\mc M_t^N(G) = \pi_t^{N,\lambda}(G) - \pi_0^{N,\lambda}(G) - \int_0^t \frac{1}{N^d} \sum_{x \in \bb T_N^d} g\big(\eta_s(x)\big) \nabla \cdot \mc
A \nabla G(x) ds
\end{equation}
plus a rest vanishing in $\mathcal{L}^{2}(\mathbb{P}_{\mu_{N}})$ as $N \to
\infty$. The next result will allow us to write the integral term
(\ref{ec3}) as a function of $\pi_t^N$ plus a vanishing term as $N \to \infty$.

\begin{proposition}(Replacement Lemma)

\quad\

 For every $\delta > 0$,
\[
\limsup_{\varepsilon \to 0} \limsup_{N \to \infty} \bb P_{\mu_N} \Big[ \int_0^T \frac{1}{N^d} \sum_{x \in \bb T_N^d} V_{\varepsilon N}
(\eta_s,x) ds > \delta\Big] = 0,
\]
where
\[
V_l(\eta,x) =\Big| \frac{1}{(2l+1)^d} \sum_{|y|\leq l} g\big(\eta(x+y)\big) - \phi\big(\eta^l(x)\big) \Big|
\]
and
\[
\eta^l(x) = \frac{1}{(2l+1)^d} \sum_{|y|\leq l} \eta(x+y).
\]
\end{proposition}

The proof of this proposition is the same as the one presented in Chapter 5 of \cite{kl}, so we omit it. Using this proposition, we see that for
any continuous function $G: U^d \to \bb R$,
\[
\int_0^t \frac{1}{N^d} \sum_{x \in \bb T_N^d} \Big\{g\big(\eta_s(x)\big) - \phi\big(\eta^{\varepsilon N}_s(x)\big)\Big\} G(x) ds \to 0
\]
in $\mathbb{P}_{\mu_{N}}$-probability as $N \to \infty$ and then $\varepsilon
\to 0$. On the other hand, since $\eta^{\varepsilon N}_s(x) =
\pi_s^N(\mathbf{1}(|u-x|\leq \varepsilon))$, we conclude that
\[
\int_0^t \frac{1}{N^d} \sum_{x \in \bb T_N^d} g\big(\eta_s(x)\big)  G(x) ds \to \int_0^t ds \int \phi(\pi(s,u))G(u)du
\]
in $\mathbb{P}_{\mu_{N}}$-probability. Since $\mc M_t^N(G)$ converges to 0, taking $N \to \infty$ in equation (\ref{ec4}) we obtain that
\begin{equation}
\label{ec5} 0 = \int\pi(t,u)G(u)du -\int \rho_0(u)G(u)du -\int_0^t\int
\phi(\pi(s,u)) \nabla \cdot \mc A \nabla G(u) du ds
\end{equation}
for every $G$ regular. Approximating a twice-differentiable function $G$ by regular functions $G_n$ in the uniform topology, we extend this
identity to functions $G \in \mc C^2(U^d)$.

Let $G:[0,T] \times U^d \to \bb R$ be of class $\mc C^{1,2}$. Take the partition $\{t_i=Ti/n; i=0,...,n\}$ of the interval $[0,T]$ and define $G_n:
[0,T] \times U^d \to \bb R$ by
\[
G_n(t,u) =\frac{n(t_i-t_{i-1})}{T} G(t_{i-1},u) + \frac{n(t_i -t)}{T} G(t_i,u),
\]
for $t \in [t_{i-1},t_i]$. In general, for a piecewise-differentiable path $G:
[0,T] \to \mc L^2_N$,
\[
\pi_t^N(G_t) - \pi_t^N(G_0) - \int_0^t\Big\{\pi_s^N(\partial_s G_s) + \frac{1}{N^d} \sum_{x \in \bb T_N^d} g\big(\eta_s(x)\big) L_N
G_s(x)\Big\}ds
\]
is a martingale of quadratic variation
\[
\int_0^t \frac{1}{N^{2d}} \sum_{x, y \in \bb T_N^d} p_N(x,y) g\big(\eta_s(x)\big) \big[G_s(y) - G_s(x)\big]^2ds.
\]

Repeating the arguments in the proof of equation (\ref{ec5}) for $G_n$, we conclude that
\[
\begin{split}
0= \int \pi(t,u) G_n(t,u)& du -\int \rho_0(u) G_n(0,u) du  \\
& - \int_0^t \int \big\{ \pi(s,u) \partial_s G_n(s,u) + \phi(\pi(s,u))\nabla
\cdot \mc A \nabla G_n(s,u) du ds.
\end{split}
\]

Taking the limit as $n$ goes to $\infty$, we obtain that
\[
\begin{split}
0= \int \pi(t,u) G(t,u) du &- \int \rho_0(u) G(0,u) du - \\
& - \int_0^t \int \big\{ \pi(s,u) \partial_s G(s,u) + \phi(\pi(s,u))\nabla \cdot
\mc A \nabla G(s,u) du ds
\end{split}
\]
for every $G:[0,T] \times U^d \to \bb R$ of class $\mc C^{1,2}$. This is the weak form of the hydrodynamic equation (\ref{echid}), see
(\ref{definitionweaksolution}). Since equation (\ref{echid}) has at most one weak solution, we conclude that $\pi(t,u) = \rho(t,u)$ $Q-a.s.$,
which ends the proof of Theorem \ref{t2}.

\section{Proof of Theorem \ref{tfeq}}
\label{s3}

Denote by $\mc Q^N$ the distribution in $\mc D([0,T], \mc S'(U^d))$ induced by the process $\mc Y_t^N$ and $\nu_{\rho}$. The standard proof of
equilibrium fluctuations cannot be accomplished for the density field $\mc Y^{N}_\cdot$. In order to overcome this problem we introduce as
before, the corrected density fluctuation field defined on functions $G\in \mc S(U^d)$ by
\begin{equation*} \label{corrdensfield}
\mc Y_{t}^{N,\lambda}(G)=\frac{1}{N^{d/2}} \sum_{x\in{\mathbb{T}_{N^d}}} G^{\lambda}_{N}(x)(\eta_{t}(x)-\rho),
\end{equation*}
where $G_{N}^{\lambda}$ is the solution of equation (\ref{poisson1}).

For $t \geq 0$, let $\mathcal{F}_{t}$ be the $\sigma$-algebra on $\mc D([0,T],\mc S'(U^d))$ generated by $\mc Y_{s}(H)$ for $s\leq{t}$ and $H$
in $\mc S(U^d)$ and set $\mathcal{F}=\sigma(\bigcup_{t\geq{0}}\mathcal{F}_{t})$. Denote by $\mc Q^{\lambda}_{N}$ the distribution on $\mc
D([0,T],\mc S'(U^d))$ induced by the corrected density fluctuation field $\mc Y_{.}^{N,\lambda}$ and $\nu_{\rho}$.

We make use of the following result, which permits to identify the limiting
process:

\begin{proposition} \label{orn-uhl}

\quad

There exists a unique process $\mc Y_{t}$ in $\mc C([0,T],\mc S'(U^d))$ such that:
\begin{itemize}
\item[i)] For every function $G\in \mc S(U^d)$,
\begin{equation*}
M_{t}(G)=\mc Y_{t}(G)-\mc Y_{0}(G)-\int^{t}_{0}\mc Y_{s}\big(\phi'(\rho)\nabla
\cdot \mathcal{A}\nabla G\big)ds
\end{equation*}
and
\begin{equation*}
(M_{t}(G))^{2}-\phi(\rho)t\int_{U^d}\nabla G(u) \cdot \mc A \nabla G(u)du
\end{equation*}
are $\mathcal{F}_{t}$-martingales.
\item[ii)]
$\mc Y_{0}$ is a Gaussian field of mean zero and covariance given by
\begin{equation} \label{cov}
E\big[\mc Y_{0}(G)\mc Y_{0}(H)\big]=\chi(\rho)\int_{U^d}G(u)H(u)du,
\end{equation}
\end{itemize}
where $\chi(\rho)=\textbf{Var}(\eta(0),\nu_{\rho})$ and $G$, $H\in{\mathcal{S}(U^d)}$. The process $\mc Y_t$ is called the generalized
Ornstein-Ulenbeck process of mean zero and characteristics $\phi'(\rho)\nabla \cdot \mathcal{A} \nabla$, $\sqrt{\phi(\rho)\mathcal{A}}\nabla$.
\end{proposition}

Theorem \ref{tfeq} is a consequence of the following result about the corrected fluctuation field.

\begin{theorem}
Let $\mc Q$ be the probability measure on $\mc C([0,T],\mc S'(U^d))$ corresponding to the stationary generalized Ornstein-Uhlenbeck process of
mean zero and characteristics $\phi'(\rho)\nabla \cdot \mathcal{A} \nabla$, $\sqrt{\phi(\rho)\mathcal{A}}\nabla$. Then the sequence $\{\mc
Q_{N}^{\lambda}\}_{N\geq{1}}$ converges weakly to the probability measure $\mc Q$. \label{th:flu}
\end{theorem}

Before we enter into the proof of this theorem, we prove Theorem
\ref{tfeq} from it. In fact, it is enough to show that
\begin{equation}
\label{ec10}
\lim_{N \to \infty} \bb E_{\nu_\rho}\big[\big(\mc Y_t^N(G) - \mc
Y_t^{N,\lambda}(G)\big)^2\big] =0
\end{equation}
for any $t \in [0,T]$, $G \in \mc S(U^d)$. But this is immediate from the fact that $G_N^\lambda$ converges to $G$ in $\mc L_N^2$ and the
independence of $\eta(x)$, $\eta(y)$ for $x\neq y$ under the invariant measure $\nu_\rho$.

In order to prove Theorem \ref{th:flu}, we need to verify that the sequence of probability measures $\{\mc Q_{N}^{\lambda}\}_{N\geq{1}}$ is
tight and to characterize the limit field. Then we show that the limit field is equal in distribution to $\mc Y_t$ using its characterization in
terms of the martingale problem (Proposition \ref{orn-uhl}).

Fix a smooth function $G \in{\mc S(U^d)}$. By Dynkin's formula,
\begin{equation}
\label{ec8} M^{N,\lambda}_{t}(G)=\mc Y^{N,\lambda}_{t}(G)-\mc Y^{N,\lambda}_{0}(G) -\int^{t}_{0} \frac{1}{N^{d/2}} \sum_{x\in\mathbb{T}_{N}^d}
g\big(\eta_{s}(x)\big) L_{N}G_{N}^{\lambda}(x)ds
\end{equation}
is a martingale with respect to the natural filtration $\mathcal F_t =\sigma(\eta_s, s\leq t)$ whose quadratic variation is given by

\begin{equation*}
\label{eq:martingale N} \<M^{N,\lambda}_{t}(G)\>=\int^{t}_{0}\frac{1}{N^d}\sum_{x,y\in\mathbb{T}_{N}^d} g\big(\eta_{s}(x)\big) p_{N}(x,y)
\big[G_{N}^{\lambda}(y) -G_{N}^{\lambda}(x)\big]^{2}ds.
\end{equation*}

At first, we establish the limit of the quadratic variation. Notice that in the previous formula we can replace $g\big(\eta_s(x)\big)$ by
$\phi(\rho)$, since
\[
\begin{split}
\bb E_{\nu_{\rho}} \Big[ \Big( \int_0^t
\frac{1}{N^{d}}
    &\sum_{x,y \in \mathbb{T}_{N}^d}
\big\{g\big(\eta_{s}(x)\big)-\phi(\rho)\big\} p_{N}(x,y)
\big[G_{N}^{\lambda}(y)-G_{N}^{\lambda}(x)\big]^{2}ds\Big)^2\Big] \leq \\
    &\leq \frac{t^{2}}{N^{2d}} \textbf{Var}(g,\nu_{\rho}) \sum_{x,y\in{\bb
    T_{N}^d}} p_{N}(x,y) \big[G_{N}^{\lambda}(y)-G_{N}^{\lambda}(x)\big]^{2}
    \times \\
    &\qquad \times \sup_{x\in{\bb T_{N}^d}}\sum_{y\in{\bb T_{N}^d}} p_{N}(x,y)
    \big[G_{N}^{\lambda}(y)-G_{N}^{\lambda}(x)\big]^{2}\\
    & \leq \frac{Ct^{2}}{N^{d-2}}
    ||R^{\lambda}G||^{2}_{0,N}||R^\lambda||_{\infty,N}.
\end{split}
\]

For dimension $d \geq 3$, this last expression goes to 0 as $N \to \infty$. In order to cover the case $d=2$, we can use Theorem (1.31) of
\cite{sz}, expression (1.32) with $t=s$ and $\alpha =1/N$ and take the Laplace transform of equation (1.32), to obtain a sharper estimate for
$G_N^\lambda(x) - G_N^\lambda(y)$. In this case, we obtain that the last line is bounded by $N^{-(d-2+2\sigma)}$, for some $\sigma>0$. As a
consequence, for any $d\geq{2}$, the quadratic variation can be written as
\begin{equation*}
\int^{t}_{0}\frac{1}{N^d}\sum_{x\in\mathbb{T}_{N}^d} \phi(\rho)\sum_{y \in \mathbb{T}_{N}^d} p_{N}(x,y)
\big[G_{N}^{\lambda}(y)-G_{N}^{\lambda}(x)\big]^{2}ds,
\end{equation*}
plus a vanishing term in the $\mc L^2(\bb P_{\nu_\rho})$-norm. Using the convergence of $G_N^\lambda$ in $\mc L^2_N$ and the resolvent estimates
in the proof of Lemma \ref{l2}, this last integral converges to
\begin{equation*}
t\phi(\rho)\int_{U}\nabla G(u) \cdot \mc A \nabla G(u) du,
\end{equation*}
as $N$ goes to $\infty$.

Now we study the limit of the martingale $M^{N,\lambda}_{t}(G)$, see expression (\ref{ec8}). Since
$\sum_{x\in{\mathbb{T}_{N}^d}}L_{N}G_{N}^{\lambda}(x)=0$, we can rewrite the integral part of the martingale as
\begin{equation*}
\int^{t}_{0}\frac{1}{N^{d/2}}\sum_{x\in\mathbb{T}_{N}^d} \{g(\eta_{s}(x))-\phi(\rho)\}L_{N}G_{N}^{\lambda}(x)ds.
\end{equation*}

On the other hand, since $G_{N}^{\lambda}$ is the solution of equation
(\ref{poisson1}), the last integral can be written as
\begin{equation*}
\int^{t}_{0} \frac{1}{N^{d/2}} \sum_{x\in\mathbb{T}_{N}^d} \{g\big(\eta_{s}(x)\big) -\phi(\rho)\} \Big\{\lambda G_{N}^{\lambda}(x)-\lambda
G(x)+\nabla\cdot\mathcal{A}\nabla G(x)\Big\}ds.
\end{equation*}

Our aim now consists in showing that it is possible to write the integral part of the martingale as the integral of a function of the density
fluctuation field plus a term that goes to zero in $\mathcal{L}^{2}(\mathbb{P}_{\nu_{\rho}})$. The first result needed to proceed in that
direction is the following:

\[
\begin{split}
\bb E_{\nu_{\rho}}\Big[ \Big( \int_0^t \frac{1}{N^{d/2}}\sum_{x\in{\mathbb{T}_{N}^d}} \{g(\eta_{s}(x))-\phi(\rho)\}
    &\big[G_{N}^{\lambda}(x)-G(x)\big]ds\Big)^2\Big] \leq \\
    & \leq C \textbf{Var}(g,\nu_\rho) ||G_N^\lambda -G||^2_{0,N}
    \xrightarrow{N \to \infty}0.
\end{split}
\]

The second one is known as the Boltzmann-Gibbs principle. Here we have the need to introduce some definitions. Take a function $f: \chi \times
\Omega_N^d \to \bb R$. For each $\omega \in \chi$ and each $x \in \bb T_N^d$, define
\[
f(x,\eta)=f(x, \eta,\omega) =: f(\theta_{xN}\omega, \tau_x\eta),
\]
where $\tau_x\eta$ is the shift of $\eta$ to $x$: $\tau_x \eta(y) = \eta(x+y)$. Notice that we do not include explicitly the dependence of
$f(x,\eta)$ in $\omega$, since in our setting $\omega$ is fixed.

\begin{definition}
We say that $f$ is local if there exists $R>0$ such that $f(\omega,\eta)$ depends only on the values of $\eta(y)$ for $|y| \leq R$. In
this case, we can consider $f$ as defined in all the spaces $\chi \times \Omega_N^d$ for $N \geq R$.
\end{definition}
\begin{definition}
 We say that $f$ is Lipschitz if there exists $c=c(\omega) >
0$ such that for all $x$, $|f(\omega, \eta)-f(\omega,\eta')| \leq c|\eta(x)-\eta'(x)|$ for any $\eta$, $\eta'$ such that $\eta(y) =\eta'(y)$ for
any $y \neq x$. If the constant $c$ can be chosen independently of $\omega$, we say that $f$ is uniformly Lipschitz.
\end{definition}

\begin{theorem}{(Boltzmann-Gibbs principle)}
\label{th:bg}

For every $G\in \mc S(U^d)$, every $t>0$ and every local, uniformly Lipschitz function $f: \chi \times \Omega_N^d \to \bb R$,
\begin{equation} \label{expBG}
\lim_{N\rightarrow{\infty}}\mathbb{E}_{\nu_{\rho}}\Big[\int_{0}^{t} \frac{1}{N^{d/2}} \sum_{x\in{\mathbb{T}_{N}^{d}}} G(x)V_{f}(x,\eta_s)ds
\Big]^{2}=0
\end{equation}
where
\begin{equation*}
V_{f}(x,\eta) =f(x,\eta)-E_{\nu_{\rho}}\big[f(x,\eta)\big]-\partial_{\rho}
E\Big[\int
f(x,\eta)d\nu_{\rho}(\eta)\Big]\big(\eta(x)-\rho\big).
\end{equation*}
Here $E$ denotes the expectation with respect to $P$, the random environment.
\end{theorem}

In order to simplify the exposition, the proof of this last result is postponed to the next section. As we need to write the integral part of
the martingale $M_{t}^{N,\lambda}(G)$ in terms of the density fluctuation field, by using the first result stated above we are able to write the
integral part of the martingale as
\begin{equation*}
\int^{t}_{0} \frac{1}{N^{d/2}} \sum_{x\in\mathbb{T}_{N}^d} \big\{g\big(\eta_{s}(x)\big)-\phi(\rho)\big\} \nabla\cdot\mathcal{A}\nabla G(x)ds
\end{equation*}
plus a term that converges to $0$ in the $\mathcal{L}^{2}(\mathbb{P}_{\nu_\rho})$-norm. The replacement of the function $g(\eta_{s})-\phi(\rho)$
by $\phi'(\rho)[\eta_{s}(x)-\alpha]$ in the last integral, is possible thanks to the Boltzmann-Gibbs principle. Doing so, the integral part of
the martingale can be written as
\begin{equation*}
M^{N,\lambda}_{t}(G)=\mc Y^{N,\lambda}_{t}(G)-\mc Y^{N,\lambda}_{0}(G) -\int^{t}_{0}\frac{1}{N^{d/2}} \sum_{x\in\mathbb{T}_{N}^d} \phi'(\rho)
\nabla\cdot\mathcal{A}\nabla G(x)\big(\eta_{s}(x)-\rho\big)ds
\end{equation*}
plus a term that vanishes in $\mc L^2(\bb P_{\nu_\rho})$ as $N \to \infty$. Notice that the integrand in the previous expression is a function
of the density fluctuation field $\mc Y_t^N$, see (\ref{densityfield}). By (\ref{ec10}), we can replace inside the integral of last expression
the density fluctuation field $\mc Y_t^N$ by the corrected density fluctuation field $\mc Y_t^{N,\lambda}$.

Suppose that the sequence $\{\mathcal{Q}_{N}^{\lambda}\}_{N\geq{1}}$ is tight and let $\mc Q^\lambda$ be a limit point of it. Denote by $\mc
Y_t$ the process in $\mc D([0,T], \mc S'(U^d))$ induced by the canonical projections under $\mc Q^\lambda$. Taking the limit as $N \to \infty$
under an appropriate subsequence in expression (\ref{ec8}), we obtain that

\begin{equation*} \label{eq:martingale M}
M^{\lambda}_{t}(G)=\mathcal{Y}_{t}(G)-\mathcal{Y}_{0}(G) -\int^{t}_{0}\mathcal{Y}_{s}(\phi'(\rho)\nabla \cdot \mathcal{A}\nabla G)ds
\end{equation*}
is a martingale of quadratic variation
\[
t \phi(\rho) \int_{U^d} \nabla G(u) \cdot \mc A \nabla G(u) du.
\]

On the other hand, it is not hard to show that $\mathcal{Y}_{0}$ is a Gaussian field with covariance given by (\ref{cov}). Therefore, $\mc
Q^\lambda$ is equal to the probability distribution $\mc Q$ of a generalized Ornstein-Uhlenbeck process in $\mc C([0,T],\mc S'(U^d))$ (and it
does not depend on $\lambda$). As a consequence, the sequence $\{ \mc Q_N^\lambda\}_{N\geq{1}}$ has at most one limit point and Theorem
\ref{th:flu} shall follow if we prove tightness for $\{\mc Q_N^\lambda\}_{N\geq{1}}$.

Lastly, it remains to treat the problem of tightness of the sequence $\{Q_{N}^{\lambda}\}_{N\geq{1}}$. For that we use a criterion due to Mitoma \cite{m} (see also \cite{fps}), which allows to conclude that the sequence is tight and that any weak limit is supported in $C([0,T],S'(U^d))$, since the following
estimates hold:

\begin{itemize}
\item[a)] For every $T>0$ and $G\in{\mc S(U^d)}$,
\begin{equation*}
\sup_{N}\sup_{t\in{[0,T]}} \bb{E}_{\nu_{\rho}}\Big[\mc Y_{t}^{N,\lambda}(G)\Big]^{2}
<\infty.
\end{equation*}

\begin{equation*}
\sup_{N} \sup_{t\in{[0,T]}} \bb{E}_{\nu_{\rho}} \Big[\frac{1}{N^{d/2}}\sum_{x\in{\bb T_{N}^d}}\phi'(\rho)\nabla \cdot \mathcal{A}\nabla
G(x)\big(\eta_{s}(x)-\rho\big)\Big]^{2}<\infty.
\end{equation*}

\begin{equation*}
\sup_{N}\sup_{t\in{[0,T]}} \bb{E}_{\nu_{\rho}} \Big[\frac{1}{N^{d/2}}\sum_{x\in{\bb T_{N}^d}}\phi(\rho) \sum_{y\in{\bb {T}_{N}^d}}
p_{N}(x,y)\big[G^{\lambda}_{N}(y)-G_{N}^{\lambda}(x)\big]^{2}\Big]^{2}<\infty.
\end{equation*}

\item[b)] For every $G\in{\mc S(U^d)}$ there exists $\delta(t,G,N)$ such that
$\lim_{N\to \infty}\delta(t,G,N)=0$ and
\begin{equation*}
\lim_{N}\bb{P}_{\nu_{\rho}}\Big(\sup_{0\leq{s\leq{t}}}\Big|\mc Y_{s}^{N,\lambda}(G)- \mc Y_{s^{-}}^{N,\lambda}(G)\Big|>\delta(t,G,N)\Big) = 0.
\end{equation*}
\end{itemize}

The first expectation in a) is bounded by $||G^{\lambda}_{N}||_{\infty}\chi(\rho)$ , which in turn is bounded by $C ||R^{\lambda}G||_{\infty}$. The
second expectation in a) is bounded by $C ||\nabla \cdot \mathcal{A}\nabla G||_{2}^{2}$ and the last one bounded by $C||R^{\lambda}G||_{2}^{4}$.

To prove b) we only have to remark that by definition of the process it holds that
$\sup_{0\leq{s\leq{t}}}|\mc Y_{s}^{N,\lambda}(G)-\mc Y_{s^{-}}^{N,\lambda}(G)\Big|\leq{\frac{||G_{N}^{\lambda}||_{\infty}}{N^{d/2}}}$.

By the results proved $\{\mc Q_{N}^{\lambda}\}_{n\geq{1}}$ is tight and we have identified above a unique limit point $\mc Q$ that corresponds
to the Ornstein-Uhlenbeck process; consequently the whole sequence converges to $\mc Q$.


\section{Boltzmann-Gibbs Principle}
\label{s4}

This section is devoted to the proof of Theorem (\ref{th:bg}). Let $f:\chi\times\Omega_{N}^d\rightarrow{+\infty}$ be a local, uniformly
Lipschitz
 function and take  $f(x,\eta)=f(\theta_{Nx}\omega,\tau_x\eta)$.

Fix a function $G\in{\mc S(U^d)}$ and an integer $K$ that shall increase to $\infty$ after $N$. For each $N$, we subdivide $\mathbb{T}_N^d$ in
non overlapping cubes of linear size $K$. Denote them by $\{I_{j}, 1\leq{j}\leq{M}^{d}\}$, where $M=[\frac{2N}{K}]$. Let $I_{0}$ be the set of
points that are not included in any $I_{j}$ which implies that $|I_{j}|\leq{dKN^{d-1}}$. If we restrict the sum in the expression that appears
inside the integral in (\ref{expBG}) to the set $I_{0}$, then its $L^{2}(\mathbb{P}_{\nu_{\rho}})$ norm clearly vanishes as
$N\rightarrow{+\infty}$.

Let $\Lambda_{s_{f}}$ be the smallest cube centered at the origin that contains the support of $f$ and define $s_f$ as the radius of
$\Lambda_{s_f}$. Denote by $I_{j}^{0}$ the interior of the interval $I_{j}$, namely the sites $x$ in $I_{j}$ that are at a distance at least
$s_{f}$ from the boundary:
\begin{equation*}
I_{j}^{0}=\{x\in{I}_{j}, d(x,\mathbb{T}_{N}^{d}\setminus{I_{j}})>{s_{f}}\}.
\end{equation*}

Denote also by $I^{c}$ the set of points that are not included in any $I_{j}^{0}$. By construction it is easy to see that
$|I^{c}|\leq{dN^{d}(\frac{c(g)}{K}+\frac{K}{N})}$.
Using the notation just settled, we have that

\begin{multline*}
\frac{1}{N^{d/2}}\sum_{x\in{\mathbb{T}_{N}^d}}
    H(x)V_{f}(x,\eta_t)=
\frac{1}{N^{d/2}}\sum_{x\in{I}^{c}}H(x)V_{f}(x,\eta_t) +\\
    +\frac{1}{N^{d/2}}\sum_{j=1}^{M^d}\sum_{x\in{I}_{j}^{0}}\Big[H(x)-H(y_{j})\Big]V_{f}(x,\eta_t)
    +\frac{1}{N^{d/2}}\sum_{j=1}^{M^d}H(y_{j})\sum_{x\in{I_{j}^{0}}}V_{f}(x,\eta_t),
\end{multline*}
where $y_{j}$ is a point in $I_{j}$. We assume that the points $y_j$ have the same relative position on each of the cubes.
The first step is to prove that
\begin{equation*}
\lim_{K\rightarrow{\infty}}\lim_{N\rightarrow{\infty}}\mathbb{E}_{\nu_{\rho}}
\Big[\int_{o}^{t}\frac{1}{N^{d/2}}\sum_{x\in{I}^{c}}H(x)V_{f}(x,\eta_t)ds\Big]^{2}=0.
\end{equation*}

Applying Schwarz inequality, since $\nu_{\rho}$ is an invariant product measure and since $V_{f}$ has mean zero with respect to the measure
$\nu_{\rho}$, the last expectation is bounded above by
\begin{equation*}
\frac{t^{2}}{N^{d}}\sum_{\substack{x,y\in{I^{c}}\\|x-y|\leq{2s_{f}}}}H(x)H(y) E_{\nu_{\rho}}\big[V_{f}(x,\eta)V_{f}(y,\eta)\big].
\end{equation*}

Since $V_{f}$ belongs to $\mathcal{L}^{2}(\nu_{\rho})$ and
$|I^{c}|\leq{dN^{d}(\frac{c(f)}{K}+\frac{K}{N})}$, the last expression
vanishes by taking first $N\rightarrow{+\infty}$ and then
$K\rightarrow{+\infty}$.

Applying the same arguments, it is not hard to show that
\begin{equation*}
\lim_{N\rightarrow{\infty}}\mathbb{E}_{\nu_{\rho}}\Big[\int_{0}^{t}\frac{1}{N^{d/2}}\sum_{j=1}^{M^d}\sum_{x\in{I}_{j}^{0}}
\big[H(x)-H(y_{j})\big]V_{f}(x,\eta_t)ds\Big]^{2}=0.
\end{equation*}

In order to finish the proof it remains to show that
\begin{equation*}
\lim_{K\rightarrow{\infty}}\lim_{N\rightarrow{\infty}} \mathbb{E}_{\nu_{\rho}}\Big[\int_{0}^{t}\frac{1}{N^{d/2}}\sum_{j=1}^{M^d}H(y_{j})
\sum_{x\in{I_{j}^{0}}}V_{f}(x,\eta_t)ds\Big]^{2}=0.
\end{equation*}

Let $\bb L_{N}$ be the generator of the zero-range process without the random environment (that is, taking $a(\omega) \equiv 1$ in (\ref{RE})),
and without the diffusive scaling $N^2$. For each $j=0,..,M^d$ denote by $\zeta_{j}$ the configuration $\{\eta(x), x \in{I_{j}}\}$ and by $\bb
L_{I_{j}}$ the restriction of the generator $\bb L_{N}$ to the interval $I_{j}$, namely:
\begin{equation*}
\bb L_{I_{j}} h(\eta) = \sum_{\substack{ x,y \in I_j \\ |x-y|= 1/N}} g\big( \eta(x) \big) \big[h(\eta^{x,y}) -h(\eta) \big].
\end{equation*}

We point out here that we are introducing a slightly different generator than the one that generates the dynamics, namely $\mathcal{L}_{N}$. The reason for
doing this stands on the fact that the dynamics generated by this operator is translation invariant. The generator that we choose to introduce
here is not random, but due to the ellipticity assumption on the environment, it is mutually bounded with the one that we have started with.

Now we introduce some notation. Fix a local function $h: \chi \times \Omega_N^d \to \bb R$, measurable with respect to
$\sigma(\eta(x),x\in{I_{1}})$, such that $E[\int h(\omega,\eta)^2 d\nu_\rho] < \infty$ and let $h_{j}$ be the translation of $h$ by $y_j -y_0$:
$h_j(x,\eta) = h(\theta_{(y_j-y_0)N} \omega, \tau_{y_j-y_0} \eta)$. Denote by $\mc L^2(\nu_\rho\times P)$ the set of such functions. Consider
\begin{equation*}
V_{H,h}^{N}(\eta)=\frac{1}{N^{d/2}}\sum_{j=1}^{M^d}H(y_{j}) \bb L_{I_{j}}h_{j}(\zeta_{j}).
\end{equation*}

By proposition A 1.6.1 of \cite{kl} and the ellipticity assumption, it is not hard to show that
\begin{equation*}
\mathbb{E}_{\nu_{\rho}}\Big[\int_{0}^{t}\frac{1}{N^{d/2}}\sum_{j=1}^{M^d}H(y_{j})
\bb L_{I_{j}}h_{j}(\zeta_{j}(s))ds\Big]^{2} \leq 20 \epsilon_0^{-1} t |||V_{H,h}^{N}|||_{-1}^2,
\end{equation*}
where the norm $|||\cdot|||_{-1}$ is given by the variational formula

\begin{equation} \label{eq1}
|||V^N_{H,h}|||_{-1}^2 = \sup_{F \in \mc L^{2}(\nu_{\rho})}\Big\{2\int V_{H,h}^{N}(\eta)F(\eta) d\nu_{\rho}-N^2\< F, -\bb L_{N}F \>_\rho \Big\},
\end{equation}
where $\<\cdot,\cdot\>_\rho$ denotes the inner product in $\mc L^2(\nu_\rho)$.

By the Cauchy-Schwarz inequality,
\begin{equation*}
\int
\bb L_{I_{j}}h_{j}(\zeta_{j})F(\eta)d\nu_{\rho} \leq {\frac{1}{2\gamma_{j}} \<-\bb L_{I_{j}}h_{j},h_{j}\>_{\rho} +\frac{\gamma_{j}}{2} \<F,-\bb L_{I_{j}}F\>_{\rho}}
\end{equation*}
for each $j$, where $\gamma_{j}$ is a positive constant. Therefore,
\begin{equation*}
2\int V^N_{H,h}(\eta) F(\eta) d \nu_\rho \leq \frac{2}{N^{d/2}}\sum_{j=1}^{M^d}H(y_{j})\Big\{\frac{1}{2\gamma_{j}} \<-\bb L_{I_{j}}h_{j},h_{j}\>_{\rho}+
\frac{\gamma_{j}}{2} \<F,-\bb L_{I_{j}}F\>_{\rho}\Big\}.
\end{equation*}

Taking for each $j$, $\gamma_{j}= N^{2+\frac{d}{2}}|H(y_{j})|^{-1}$ we have that
\begin{equation*}
\frac{2}{N^{d/2}}\sum_{j=1}^{M^d}\Big|H(y_{j})\Big|\frac{\gamma_{j}}{2} \<F,-\bb L_{I_{j}}F\>_{\rho}\leq N^{2}\<F,-\bb L_{N}F\>_{\rho},
\end{equation*}
and the expectation becomes bounded by
\begin{equation*}
\frac{20\epsilon_0^{-1} t}{N^{d/2}} \sum_{j=1}^{M^d} \frac{\big|H(y_{j})\big|}{\gamma_{j}} \<-\bb L_{I_{j}}h_{j},h_{j}\>_{\rho}
    \leq \frac{20tM^d||H||_\infty}{\epsilon_0^2 N^{2+d}} \sum_{j=1}^{M^d} \frac{1}{M^d} \<-\bb L_j h_j,h_j\>_\rho.
\end{equation*}

By the ergodic theorem, the sum in the previous expression converges as $N \to \infty$ to a finite value and therefore this last expression
vanishes as $N\rightarrow{\infty}$. To conclude the proof of the theorem we need to show that
\begin{multline*}
\lim_{K\rightarrow{\infty}}\inf_{h\in{\mc L^{2}(\nu_{\rho}\times P)}} \\
\lim_{N\rightarrow{\infty}}\mathbb{E}_{\nu_{\rho}}\Big[\int_{0}^{t}\frac{1}{N^{d/2}}
\sum_{j=1}^{M^d}H(y_{j})\Big\{\sum_{x\in{I_{j}^{0}}}V_{f}(x,\eta_s)-\bb L_{I_{j}}h_{j}(\zeta_{j}(s))\Big\}\Big]^{2}=0.
\end{multline*}

By Schwarz inequality the expectation in the previous expression is bounded by
\begin{equation*}
\frac{t^{2}}{N^{d}}\sum_{j=1}^{M^d}||H||_\infty^{2} E_{\nu_{\rho}}\Big(\sum_{x\in{I_{j}^{0}}}V_{f}(x,\eta)-\bb L_{I_{j}}h_{j}(\zeta_{j})\Big)^{2}
\end{equation*}
because the measure $\nu_{\rho}$ is invariant under the dynamics and also translation invariant and the supports of $V_{f}(x,\eta)-\bb
L_{I_{i}}h_{i}(\zeta_{i})$ and $V_{f}(y,\eta)-\bb L_{I_{j}}h_{j}(\zeta_{j})$ are disjoint for $x\in{I_{i}^{0}}$ and $y\in{I_{j}^{0}}$,
 with $i\neq{j}$.

By the ergodic theorem, as $N\rightarrow{\infty}$ this expression converges to
\begin{equation}
\frac{t^{2}}{K^{d}}||H||_\infty^2 E\Big[\int \Big(\sum_{x\in{I_{1}^{0}}}V_{f}(x,\eta)-\bb L_{I_{1}}h(\omega,\eta)\Big)^{2} d\nu_\rho\Big].
\label{eq:bg2}
\end{equation}

So it remains to show that

\begin{equation*}
\lim_{K\rightarrow{\infty}} \frac{t^{2}}{K^{d}} ||H||_\infty^2
    \inf_{h\in{\mc L^{2}(\nu_{\rho} \times P)}}
    E\Big[ \int \Big(\sum_{x\in{I_{1}^{0}}}V_{f}(x,\eta)
    -\bb L_{I_{1}}h(\omega, \eta)\Big)^{2} d \nu_\rho \Big]=0.
\end{equation*}

Denote by $\mathcal{R}(\bb L_{I_{1}})$ the range of the generator $\bb L_{I_{1}}$ in $\mc L^{2}(\nu_{\rho} \times P)$ and by
$\mathcal{R}(\bb L_{I_{1}})^{\perp}$ the space orthogonal to $\mathcal{R}(\bb L_{I_{1}})$. The infimum  of (\ref{eq:bg2}) over all
$h\in{\mc L^{2}(\nu_{\rho} \times P)}$ is equal to the projection of $\sum_{x\in{I_{1}^{0}}}V_{f}(x,\eta)$ into $\mathcal{R}(\bb L_{I_{1}})^{\perp}$.

It is not hard to show that $\mathcal{R}(\mathbb{L}_{I_{1}})^{\perp}$ is the space of functions that depends on $\eta$ only through the total
number of particles on the box $I_1$. So, the previous expression is equal to
\begin{equation}
\label{ec9} \lim_{K \to \infty} \frac{t^{2}||H||_\infty^2}{K^{d}} E\Big[\int
\Big(E_{\nu_{\rho}}\Big[\sum_{x\in{I_{1}^{0}}}V_{f}(x,\eta)\Big|\eta^{I_{1}}\Big]\Big)^{2} d\nu_\rho \Big]
\end{equation}
where $\eta^{I_{1}}=K^{-d}\sum_{x\in{I_{1}}}\eta(x)$.

Let us call this last expression $\mc I_0$. Define $\psi(x,\rho) = E_{\nu_\rho}[ f(\theta_x \omega)]$. Notice that $V_f(x,\eta)= f(x,\eta) -
\psi(x,\rho) - E[\psi'(x,\rho)]\big(\eta(x)-\rho\big)$, since in the last term the derivative commutes with the expectation with respect to the
random environment. In order to estimate the expression (\ref{ec9}) using the elementary inequality $(x+y)^2\leq{2x^2+2y^2}$, we split it into three
pieces: $\mc I_0 \leq 4(\mc I_1+ \mc I_2 +\mc I_3)$, where
\[
\mc I_1 = \frac{1}{K^d} E\Big[ \int \Big(\sum_{x \in I_1^0} E_{\nu_\rho} \big[f(x,\eta)|\eta^{I_1}\big] - \psi(x,\eta^{I_1})\Big)^2 d\nu_\rho\Big],
\]
\[
\mc I_2 =  \frac{1}{K^d} E\Big[ \int \Big(\sum_{x \in I_1^0} \psi(x,\eta^{I_1})
-\psi(x,\rho) - \psi'(x,\rho)[\eta^{I_1}-\rho] \Big)^2 d\nu_\rho \Big],
\]
\[
\mc I_3 = \frac{1}{K^d} E\Big[ E_{\nu_\rho}\Big[ \Big(\sum_{x \in I_1^0}\big(\psi'(x,\rho) - E[\psi'(x,\rho)]\big)\big[\eta^{I_1} - \rho\big] \Big)^2\Big]\Big].
\]

We will make use of the following lemma, known as the equivalence of ensembles.

\begin{lemma}
\label{eqens} Let $h: \Omega_N^d \to \bb R$ a local, uniformly Lipschitz function. Then, for each $\beta \geq 0$ there exists a constant $C$
that depends on $h$ only through its support and its Lipschitz constant, such that
\[
\Big| E_{\nu_\rho} [ h(\eta)|\eta^N] - E_{\nu_{\eta^N}}[h(\eta)]\Big| \leq \frac{C}{N^d}
\]
whenever $\rho, \eta^N \leq \beta$, where
\[
 \eta^N = \sum_{|x|\leq N} \eta(x).
\]
\end{lemma}

In order to estimate $\mc I_1$ and $\mc I_2$, we introduce the indicator functions $\mathbf{1}(\eta^{I_1} \leq \beta)$. By a large deviations
estimate, $\nu_\rho(\eta^{I_1} \geq \beta) \leq \exp(-C(\beta)K^d)$. Since $f$ is Lipschitz, it has bounded exponential moments of any order and a simple
Schwarz estimate shows that we can introduce the indicator function $\mathbf{1}(\eta^{I_1} \leq \beta)$ into the integrals in $\mc I_1$ and $\mc
I_2$. By Lemma \ref{eqens},
\[
\frac{1}{K^d} E\Big[ \int \Big(\sum_{x \in I_1^0} E_{\nu_\rho} \big[f(x,\eta)|\eta^{I_1}\big] -
\psi(x,\eta^{I_1})\Big)^2 \mathbf{1}(\eta^{I_1} \leq \beta) d\nu_\rho\Big] \leq \frac{C}{K^d},
\]
which vanishes as $K \to \infty$.

 Using a Taylor expansion for
$\psi(x,\rho)$, we see that
\[
\frac{1}{K^d} E\Big[ \int \Big(\sum_{x \in I_1^0} \psi(x,\eta^{I_1})
-\psi(x,\rho) - \psi'(x,\rho)[\eta^{I_1}-\rho] \Big)^2 d\nu_\rho
\Big] \leq \frac{C}{K^d}
\]
and also goes to 0 as $K \to \infty$.

Finally, we see that
\[
\mc I_3 = E_{\nu_\rho}\big[(\eta(0) -\rho)^2\big]\cdot E\Big[\Big( \frac{1}{K^d} \sum_x (\psi'(x,\rho) - E[\psi'(x,\rho)]\Big)^2\Big]
\]
and it goes to 0 as $K \to \infty$ by the $\mc L^2$-ergodic theorem.

\subsection{Some applications of the Boltzmann-Gibbs principle}

In the proof of Theorem \ref{tfeq}, we need to use the Boltzmann-Gibbs Principle \ref{th:bg} for the function $g(\eta(0))$, that does not
depend on the random environment. In particular, the results of the previous section are not needed in the proof of Theorem \ref{tfeq},
since the proof for the non-random case applies directly for functions that do not depend on the random environment. We point out here two
applications for the Boltzmann-Gibbs principle as stated in Theorem \ref{th:bg}.

\textbf{First application:} Consider, for simplicity, some local, bounded and uniformly Lipschitz function $f(\omega,\eta)$ that does not depend on the value of $\eta(0)$. For each $\eta \in \Omega_N$, define
\[
\Theta_x^+ \eta(z) = \begin{cases}
                      \eta(x)+1, z=x\\
              \eta(z), z \neq x,
                     \end{cases}
\]
\[
\Theta_x^- \eta(z) = \begin{cases}
                      \eta(x)-1, z=x\\
              \eta(z), z \neq x.
                     \end{cases}
\]
Notice that $\Theta_x^- \eta$ is well defined only if $\eta(x) \geq 1$. We can define a reaction-diffusion model adding to the zero-range dynamics a Glauber dynamics as follows:
\[
\begin{split}
\mc L_N^{rd}F(\eta) =: N^2 \mc L_N &+ \sum_{x \in \bb T_N^d} f(x,\eta) \big[F(\Theta_x^+ \eta) - F(\eta)\big] \\
    &+ \sum_{x \in \bb T_N^d} \alpha(\rho) \frac{f(x,\Theta_x^-\eta)}{g(\eta(x))} \big[F(\Theta_x^- \eta) - F(\eta)\big],\\
\end{split}
\]
where we define $f(x,\Theta_x^-\eta)/g(\eta(x))=0$ if $\eta(x)=0$. We have chosen the annihilation rate in such a way that the measure $\nu_\rho$ is invariant for this process. Therefore, we can obtain the equilibrium fluctuations for this model as in Section \ref{s3}.

\textbf{Second application:} This one has to do with the convergence of additive functionals of Markov processes. For each $f$ satisfying the
conditions of Theorem \ref{th:bg}, define the density fluctuation field for $f$ acting on functions $G\in \mc S(U^d)$ as
\begin{equation*}
\mathcal{Z}_{t}^{N,f}(G)=\frac{1}{N^{d/2}}\sum_{x\in{\mathbb{T}_{N}^d}}G(x))\big\{f(x,\eta_s) - E_{\nu_\rho} [f(x,\eta)]\big\}.
\end{equation*}
Note that for $f(x,\eta)=\eta(x)-\rho$, the density fluctuation field for $f$ is the density fluctuation field introduced above and denoted by
$\mathcal{Y}^{N}_s(G)$.

For fixed $f$ as above, define the additive functional
\[
 \mc I_f^N(t) = \int_0^t \mathcal{Z}_{t}^{N}(G) ds.
\]

Then, by Theorems \ref{tfeq} and \ref{th:bg},
\[
\lim_{N \to \infty} \mc I^N_f(t)  = \partial_\rho E\Big[ \int f(\omega,\eta)d \nu_\rho\Big] \int_0^t \mc Y_s(G)ds \text{ in distribution. }
\]

\begin{appendix}
\section{Some estimates for $\eta_t$}

\subsection{Entropy production}

Denote by $\mu_{N}(t)=S_t^N \mu_N$ the distribution of $\eta_t$ in $\Omega_N^d$ under $\bb P_{\mu_N}$ and define $f_t^N = \frac{d\mu_{N}(t)}{d
\nu_{\bar \rho}}$. The density $f_t^N$ satisfies the Kolmogorov equation
\[
\frac{d}{dt} f_t^N(\eta) = \mc L_N f_t^N(\eta).
\]

For each density $f: \Omega_N^d \to \bb R_+$, define the Dirichlet form $\mc D_N(f)$ by
\[
\mc D_N(f) = \sum_{\substack{x, y \in \bb T_N\\|x-y|=1/N}} \int g\big(\eta(x)\big)\big[\sqrt{f(\eta^{xy})}- \sqrt{f(\eta)}\big]^2 d\nu_{\bar \rho},
\]
and the entropy $H_N(f) = \int f\log f d\nu_{\bar \rho}$.
By the ellipticity assumption in $p_N(x,y)$, the entropy production is bounded by the Dirichlet form of $f_t^N$ \cite{kl}:
\[
\frac{d}{dt}H_N(f_t^N) \leq - 2\epsilon_0 N^2 \mc D_N(f_t^N).
\]
Assume that $H_N(\mu_N|\nu_{\bar \rho}) \leq K_0 N^d$, or in other words that $H_N(f_0^N) \leq K_0 N^d$. Since the Dirichlet form and
the entropy are convex functions of $f$, integrating the previous inequality we obtain the bounds
\[
H_N(\bar f_T^N) \leq \frac{K_0}{T}N^d, \quad\   \mc D_N(\bar f_T^N) \leq \frac{K_0}{2 \epsilon_0 T} N^{d-2},
\]
where
\[
\bar f_T^N(\eta) = \frac{1}{T} \int_0^T f_t^N(\eta) dt.
\]

\subsection{Attractiveness of $\eta_t$}

Take two probability measures $\mu$, $\nu$ in $\Omega_N$ such that $\mu \preceq \nu$. When the jump rate $g(\cdot)$ is non-decreasing, it is possible to construct a process $(\eta_t,\eta_t')$ in $\Omega_N \times \Omega_N$, starting from a coupling $\bar \mu$ of $\mu$ and $\nu$, such that for every $t \in [0,T]$
\begin{itemize}
\item[i)] The distribution of $\eta_t$ in $\mc D([0,T],\Omega_N)$ is equal to $\bb P_\mu$.
\item[ii)] The distribution of $\eta'_t$ in $\mc D([0,T],\Omega_N)$ is equal to $\bb P_\nu$.
\item[iii)] The distribution of $(\eta_t,\eta'_t)$ in $\mc D([0,T],\Omega_N \times \Omega_N)$ is concentrated on the set $\{(\eta,\eta') \in \Omega_N \times \Omega_N; \eta \preceq \eta'\}$.
\end{itemize}

In this case the process $\eta_t$ is said to be attractive. We say that a function $h: \Omega_N^d \to \bb R$ is non-decreasing if for $\eta \preceq \eta'$
then $h(\eta) \leq h(\eta')$. The following proposition is an immediate consequence of the existence of the process $(\eta_t,\eta_t')$.

\begin{proposition}
\label{p2} Let $\mu$, $\nu$ be two probability measures in $\Omega_N^d$ such that $\mu \preceq \nu$. Let $h: \Omega_N^d\rightarrow{\mathbb{R}}$
be a non-decreasing function. Then,

\[
\bb E_{\mu}\big[h(\eta_t)\big] \leq \bb E_{\nu}\big[h(\eta_t)\big]
\]
for all $t \in [0,T]$.
\end{proposition}

\subsection{An $\mc L^2$ estimate for $\eta_t$}

Consider the process $\eta_t$ starting from the equilibrium measure $\nu_\rho$. Define the $\mc L^2_N$-norm of $\eta_t$ by
\[
||\eta_t||_{0,N}^2 = \frac{1}{N^d} \sum_{x \in \bb T_N} \eta_t(x)^2.
\]

By Dynkin's formula,
\begin{equation}
\label{ec6} \mc M_t^N = ||\eta_t||_{0,N}^2 -||\eta_0||_{0,N}^2 - \int_0^t \mc L_N||\eta_s||_{0,N}^2 ds
\end{equation}
is a martingale of quadratic variation
\[
\<\mc M_t^N\> = \int_0^t ||\eta_s||_{0,N}^2(-\mc L_N)||\eta_s||_{0,N}^2 ds.
\]

Explicit computations show that $\bb E_{\nu_\rho}\big[\<\mc M_t^N\>\big] \leq C/N^{d-2}$. Therefore, by Doob's inequality,
\[
\bb E_{\nu_\rho}\Big[ \sup_{t \in [0,T]} |\mc M_t^N|^2\Big] \leq C/N^{d-2}.
\]

For the integral term in \ref{ec6}, we have the following estimate:
\[
\bb E_{\nu_\rho}\Big[ \big(\sup_{t\in [0,T]} \int_0^t \mc L_N ||\eta_s||_{0,N}^2 ds\big)^2\Big] \leq Ct E_{\nu_\rho} \Big[ ||\eta||_{0,N}^2
(-\mc L_N) ||\eta||_{0,N}^2\Big].
\]

Therefore, for dimension $d \geq 2$, we conclude that $\bb E_{\nu_\rho}\big[\sup_t ||\eta_t||_{0,N}^2\big]$ is uniformly bounded in $N$. Since
$||\eta||_{0,N}^2$ is an increasing function, we have proved the following result:

\begin{lemma}
\label{l1} Fix $\rho>0$. Let $\{\mu_N\}_{N\geq{1}}$ be a sequence of measures such that $\mu_N \preceq \nu_\rho$ for all $N$. Then,
\[
\sup_{N \in \bb N} \bb E_{\mu_N} \big[ \sup_{t \in [0,T]} ||\eta_t||_{0,N}^2 \big] < +\infty.
\]
\end{lemma}

\section{Functional analysis in the spaces $\mc L^2_N$, $\mc H_{1,N}$}
\subsection{Convergence in $\mc L^2_N$, $\mc H_{1,N}$}

\quad\

 \label{ap1} Fix $f \in \mc H_{1,N}$. We define the linear interpolation $\mc T_N^1 f$ of $f$ as follows. To fix ideas, take $d=3$. We
divide each of the cubes of size $1/N$ in $\bb T_N^d$ into six tetrahedrons with vertices in $\bb T_N^d$. The way we do this is not important,
but we do it in the same way for every cube in $\bb T_N^d$.

For a point $u$ in one of such tetrahedrons, we define $\mc T_N^1 f(u)$ as the linear interpolation of the values of $f$ on the vertices of the
tetrahedron. In this way we have defined a function $\mc T_N^1 f$ in $\mc H_1(U^d)$.

We say that $f_N \in \mc H_{1,N}$ converges strongly (resp. weakly) in $\mc H_{1,N}$ to $f \in \mc H_1(U^d)$ if
\[
\lim_{N \to \infty} \mc T_N^1 f_N = f \text{ strongly (resp. weakly) in }  \mc H_1(U^d).
\]

In an analogous way, for each $u \in U^d$ we define $\mc T_N^0 f(u)= f(x)$ if $|u-x| \leq 1/2N$. We say that $f_N$ converges strongly (resp.
weakly) in $\mc L^2_N$ to $f \in \mc L^2(U^d)$ if $\mc T_N^0 f_N$ converges strongly (resp. weakly) to $f$ in $\mc L^2(U^d)$.

A sequence $f_N \in \mc H_{-1,N}$ converges  to $f \in \mc H_{-1}(U^d)$ strongly (resp. weakly) if for any sequence $g_N \in \mc H_{1,N}$ and $g
\in \mc H_1(U^d)$ such that $g_N \to g$ weakly in (resp. strongly) $\mc H_{1,N}$ we have
\[
\lim_{N \to \infty} \<f_N, g_N\>_N = \<f,g\>.
\]

\subsection{Resolvent estimates}

\quad\

Let $f$ be a regular function and let $u_N$ be the solution of the resolvent equation
\begin{equation}
\label{ec7}
\lambda u_N(x) - L_N u_N(x) = f(x).
\end{equation}

\begin{lemma}
\label{l2}
There exists a constant $c=c(\lambda)$ such that
\[
\max\{||u_N||_{0,N},||u_N||_{1,N}\} \leq c||f||_1.
\]
\end{lemma}
\begin{proof}
By Lax-Milgram's lemma, this equation has a unique solution in $\mc H_{1,N}$. Taking the inner product of equation (\ref{ec7}) with respect to
$u_N$, we see that
\[
\lambda ||u_N||_{0,N}^2 + \<u_N, -L_N u_N\>_N \leq \<f, u_N\>_N.
\]

By the Cauchy-Schwarz inequality, $|\<f,u_N\>_N \leq ||u_N||_{0,N}||f||_{0,N}$. Using the ellipticity assumption, we obtain the estimates
\[
\begin{split}
||u_N||_{0,N} &\leq \lambda^{-1} ||f||_{0,N} \\
||u_N||_{1,N}^2 &\leq \big[(\lambda \epsilon_0)^{-1} +\lambda^{-2}\big]||f||_{0,N}^2.
\end{split}
\]

By the finite elements theory \cite{b}, there exists a constant $\gamma$ independent of $N$ such that for every $f \in \mc H_1(U^d)$,
$||f||_{0,N} \leq \gamma||f||_1$. Therefore, it
 is enough to take $c= \gamma \max\{\lambda^{-1},(\lambda \epsilon_0)^{-1} +\lambda^{-2}\}$.
\end{proof}
Since the operator $L_N$ is the generator of a random walk in $\bb T_N^d$, the solutions of (\ref{ec7}) satisfy the maximum principle:
\[
\inf_{x \in \bb T_N^d} \lambda^{-1} f(x) \leq \inf_{x \in \bb T_N^d} u_N(x) \leq \inf_{x \in \bb T_N^d}  u_N(x) \leq \sup_{x \in \bb T_N^d}
\lambda^{-1}f(x).
\]

In particular, for $f$ continuous, $||u_N||_\infty \leq \lambda^{-1} ||f||_\infty$.
\end{appendix}

\end{document}